\documentclass[11pt]{amsart}
\usepackage{amsfonts}

\usepackage{graphicx}
\newtheorem{theorem}[subsection]{Theorem}
\newtheorem{lemma}[subsection]{Lemma}

\newenvironment{ulemma}{\par\noindent\textbf{Lemma}\,\,\em}{\rm}

\newenvironment{captheorem}{\par\noindent\textbf{Captures Theorem}\,\,\em}{\rm}

\newenvironment{mattheorem}{\par\noindent\textbf{Matings Theorem}\,\,\em}{\rm}

\begin{document}

\title{Multiple equivalent matings  with the aeroplane polynomial}

\author{Mary Rees}
\address{Dept. of Math. Sci., University of Liverpool, Mathematics Building, Peach St., Liverpool L69 7ZL, U.K..}
\email{maryrees@liv.ac.uk}
\def\Box{\hbox{$\sqcap \unskip \kern -6.5pt 
\sqcup$}}
\def \Amalg{\mathbin{\raise .5pt%
	\hbox{$\scriptstyle \amalg$}}}

\numberwithin{equation}{subsection}

\begin{abstract}
We produce arbitrarily large equivalence  classes of  matings with the aeroplane polynomial. These are obtained by a slight generalisation of the technique of proof of a similar result for Wittner captures.
\end{abstract}
\begin{center}
{\em Dedicated to the memory of Adrien Douady}
\end{center}

 \maketitle
 
 \section{Introduction}\label{1}
\subsection{}\label{1.1}
Mating of polynomials is a construction invented by A. Douady and J.H. Hubbard in the 1980's, soon after their seminal work on the structure of the Mandelbrot set for quadratic polynomials (\cite{D-H1, D-H2}). The construction is of a rational map $f_{1}\Amalg f_{2}$ from any two given critically finite quadratic polynomials $f_{1}$ and $f_{2}$, apart from some restrictions on $(f_{1},f_{2})$ which are (with hindsight) very natural and simple. The construction suggests that the intricacy of structure of the parameter space of quadratic polynomials is, at the very least, reproduced in the parameter space of rational maps. Given that the coordinates of a mating are polynomials, one might expect a product structure of some sort in the space of rational maps. It has been known for some time (from unpublished work of Adam Epstein) that this cannot be true in any strong topological sense. But the way in which the structure translates to the parameter space of rational maps is even more interesting than might at first appear. This has been one aspect of an extensive study. See, for example \cite{R1, R2, R3, R5, R6}. In the current work, we avoid direct reference to such theory, and present some examples. The construction of these examples is self-contained, although they were found during a larger exploration of the theory.  We show that, even if we restrict to a single choice of $f_{1}$, the so-called aeroplane polynomial, the map 
$$f_{2}\mapsto  f_{1}\Amalg f_{2}$$
is unboundedly many-to-one.  Thus, a critically finite quadratic rational map can be realised as a mating of the aeroplane polynomial with another polynomial in arbitrarily many ways. We also obtain a similar result for {\em{captures}}. The capture construction probably first appears in written form in Wittner's 1988 thesis \cite{W}, but was studied by the Hubbard group earlier in the 1980's. (The term, ``capture'' is used here in a more restrictive sense than is sometimes the case. The precise definition will be given later.)  The capture  $\sigma _{\beta _{x}}\circ f$ is a critically finite quadratic rational map made from a critically finite quadratic polynomial $f$ and a path $\beta _{x}$ from the fixed critical point $\infty $ of $f$, which crosses the Julia set of $f$ just once, at a point $x$ in the boundary of a preperiodic Fatou component, and ends at the precritical point in that Fatou component. As in the case of the mating construction, there are some extra, natural and easily stated restrictions.  One of the constant refrains of  \cite{R1, R2,R3, R5, R6} is that ``type III'' maps, of which captures are examples, are technically easier to work with than ``type IV'' maps, of which matings are examples. It is therefore no surprise that  the mating result is essentially obtained from the capture result. The way in which this is done is reminiscent of argument in \cite{R1}. It is not difficult and may indicate that more general theoretical results can be obtained in the same way. 

I should like to thank the referee for full and prompt engagement with the paper, and for some helpful comments.

\subsection{Lamination maps}\label{1.2}
We use the notation and conventions of \cite{R5}. This means, in particular, that we work with the lamination map which is Thurston equivalent to the aeroplane polynomial, rather than with the aeroplane polynomial  itself. We call this lamination map $s_{3/7}$. It preserves the unit disc, unit circle and exterior of the unit disc. On the exterior of the unit disc, $s_{3/7}(z)=z^{2}$, and on the interior of the unit disc $s_{3/7}$ preserves a lamination known as $L_{3/7}$, and also preserves the set of complementary components of $\cup L$ in the unit disc. This lamination is an example of a {\em{quadratic invariant lamination}}. It is a closed set of disjoint chords, or {\em{leaves}}, in the unit disc. {\em{Invariance}} of a lamination $L$ in the unit disc means that:
\begin{itemize}
\item whenever there is a leaf of $L$ joining $z_{1}$ and $z_{2}$, there is also a leaf of $L$ joining $z_{1}^{2}$ and $z_{2}^{2}$;
\item whenever there is a chord joining $z_{1}$ and $z_{2}$, there are points $\pm z_{3}$ and $\pm z_{4}$ with $z_{3}^{2}=z_{1}$ and $z_{4}^{2}=z_{2}$, and such that there are leaves of $L$ joining $z_{3}$ to $z_{4}$, and $-z_{3}$ to $-z_{4}$.
\end{itemize} 

The label $3/7$ on $L_{3/7}$ means that the (common) image of the  two longest leaves of $L_{3/7}$ has one endpoint at $e^{2\pi i(3/7)}$. This image is known as the {\em{minor leaf}} of $L_{3/7}$. The other endpoint of the minor leaf is at $e^{2\pi i(4/7)}$. In what follows, we only consider invariant  laminations in which every leaf  has exactly two preimages. Such an invariant lamination is uniquely determined by its minor leaf. We therefore write $L_{3/7}=L_{4/7}$ and $s_{3/7}=s_{4/7}$. The minor leaf is periodic under $s_{3/7}$, of period $3$. It is adjacent to a unique complementary component of $\cup L_{3/7}$ in the unit disc, also known as a {\em{gap}} of $\cup L_{3/7}$,  which is also of period three. The gap adjacent to the minor leaf has a connected  preimage under $s_{3/7}$. this connected preimage contains $0$ and is known as  the {\em{critical gap}}, since it maps with degree two onto its image. It then maps homeomorphically  on the rest of the cycle, and returns to itself under $s_{3/7}^{3}$ with degree two. This is consistent with $s_{3/7}$ having a unique critical point in the unit disc at $0$, and this critical point being of period $3$. We choose $s_{3/7}$ like this. 

More generally, $L_{p}$ and $s_{p}$ are uniquely determined in the same way, for any odd denominator rational $p$ in $[0,1]$. So $s_{p}$ is a critically periodic branched covering, with one fixed critical point (at $\infty $), and one periodic critical point of the same period as the period of $p$ under $x\mapsto 2x{\rm{\ mod\ }}1$. There is a one-to-one correspondence between such maps $s_{p}$ and critically periodic polynomials in the family $z^{2}+c$ (for $c\in \mathbb C$), where the correspondence is defined by {\em{Thurston equivalence}}. In this paper, two critically finite branched coverings $f_{0}$ and $f_{1}$ with numbered postcritical sets are {\em{(Thurston) equivalent}} if there is a homotopy through critically finite branched coverings $f_{t}$ such that the postcritical set $X(f_{t})$ moves isotopically and the numbering of $X(f_{0})$ and $X(f_{1})$ is preserved by this isotopy. The stipulation that numbering is preserved is not universal, but is convenient here, especially when we consider matings. For quadratic branched coverings with a fixed critical point it is superfluous, because a fixed critical point moves isotopically through fixed critical points.

The {\em{aeroplane polynomial}} is the (unique) quadratic polynomial of the form $z^{2}+c$ which is Thurston equivalent to $s_{3/7}$.The critical point $0$ is of period $3$.  The parameter value $c$ is real and in $(-2,-1)$, and is the unique such value for which $0$ is of period $3$.

\subsection{Captures and matings}\label{1.3}

For any path $\beta $, we define a homeomorphism $\sigma _{\beta }$, up to isotopy of a suitable type. The definition can be made for any path, but since we only need the definition for paths which are arcs, we only give the definition for these. This homeomorphism $\sigma _{\beta }$ is the identity outside a suitably small disc neighbourhood of $\beta $, and maps the startpoint of $\beta $ to the endpoint, as shown in Figure 1.
\begin{figure}
\centering{\includegraphics[width=4cm]{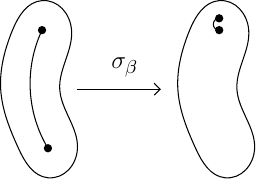}}
\caption{$\sigma _{\beta }$}
\end{figure}

A {\em{capture path for $s_{p}$}} is a path from $\infty $ to a nonperiodic point $x$ in the backward orbit of the (periodic) point $0$ which crosses the unit circle just once, into a gap of $L_{p}$ and ending at $x$, with no intersections with $L_{p}$. A {\em{capture by $s_{p}$}} is then a branched covering of the form $\sigma _{\beta }\circ s_{p}$ for a capture path $\beta $. Such a branched covering is critically finite, with critical points at $0$ and $\infty $ and critical values at $s_{p}(0)$ and $x$, the endpoint of $\beta $. Then $0$ is periodic under $\sigma _{\beta }\circ s_{p}$, and  the critical point $x$ is eventually mapped by $\sigma _{\beta }\circ s_{p}$ to the periodic orbit of $0$. We refer to such a degree two branched covering, with one critical point periodic, and the other nonperiodic, but in the backward orbit of the first, as being {\em{of type III}}.

In this paper, we only consider captures by $s_{3/7}$, which, up to equivalence, are {\em{captures by the aeroplane polynomial}}, since they can equally well be described  (up to equivalence) by capture paths in the dynamical plane of the aeroplane polynomial, that is, paths which cross the Julia set of the aeroplane polynomial just once. By the result of Tan Lei, or at least the technique of it (\cite{TL}, see also \cite{R1}) a capture $\sigma _{\beta }\circ s_{3/7}$ is Thurston equivalent to a quadratic rational map (which is unique up to M\"obius conjugation) if and only if the endpoint $x$ of the path is a point in the full orbit of $0$ which is in the larger region of the unit disc bounded by the leaf of $L_{3/7}$ with endpoints at $e^{\pm 2\pi i(1/3)}$.  This means that the crossing point by the capture path  of the unit circle is at $e^{2\pi iy}$ for $y\in (-\frac{1}{3},\frac{1}{3})$.  In section 2.8 of  \cite{R5}, elementary equivalences between different captures by $s_{3/7}$ were discussed, for capture paths with the same endpoint. It was shown that any two captures $\sigma _{\beta _{1}}\circ s_{3/7}$ and  $\sigma _{\beta _{2}}\circ s_{3/7}$, such that $\beta _{1}$ and $\beta _{2}$ have the same endpoint,  are Thurston equivalent, except when the common endpoint of the paths is in the region of the unit disc bounded between the leaves with endpoints $e^{\pm 2\pi i(1/3)}$, and endpoints $e^{\pm 2\pi i(2/7)}$. In the exceptional case, $\sigma _{\beta _{1}}\circ s_{3/7}$ and  $\sigma _{\beta _{2}}\circ s_{3/7}$ are still Thurston equivalent if $\beta _{1}$ and $\beta _{2}$ have the same endpoint, and if the $S^{1}$-crossing points $e^{2\pi iy_{1}}$ and $e^{2\pi iy_{2}}$ are such that $y_{1}$ and $y_{2}$ are both in $(\frac{2}{7},\frac{1}{3})$, or both in $(\frac{2}{3},\frac{5}{7})$. An example was given  in 9.2 of \cite{R5}  of capture paths $\beta _{1}$ and $\beta _{2}$ with  different endpoints, with crossing points $e^{2\pi iy_{1}}$ and $e^{2\pi iy_{2}}$ with $y_{1}\in (\frac{2}{7},\frac{1}{3})$ and $y_{2}\in (\frac{2}{3},\frac{5}{7})$, such that  $\sigma _{\beta _{1}}\circ s_{3/7}$ and  $\sigma _{\beta _{2}}\circ s_{3/7}$ are Thurston equivalent. It was also shown in the ``easy'' Main Theorem of \cite{R5} (2.10 and 6.1) that captures with different endpoints are not Thurston equivalent except when the endpoints are in the region bounded by the leaves joining $e^{\pm 2\pi i(1/3)}$ and $e^{\pm (2\pi i(2/7)}$. (This is essentially a folklore result if the crossing point $e^{2\pi iy}$ satisfies $y\in (\frac{1}{7},\frac{2}{7})\cup (\frac{5}{7},\frac{6}{7})$, but not so if $y\in (-\frac{1}{7},\frac{1}{7})$.) A detailed analysis of the {\em{parameter capture map}} for the aeroplane polynomial
$$[\beta ] \mapsto [\sigma _{\beta }\circ s_{3/7}]$$
--- where $[.]$ on the domain is suitable homotopy equivalence, and $[.]$ and on the range is Thurston equivalence --- will be given in \cite{R6}. The purpose of the present paper  is simply to give examples to show that the fibres of this capture map --- that is, the sets of equivalent captures --- can be arbitrarily large.

We also describe matings only up to Thurston equivalence. For (quadratic) lamination maps $s_{p}$ and $s_{q}$,  the {\em{mating}} $s_{p}\Amalg s_{q}$ of $s_{p}$ and $s_{q}$ is the branched covering defined by
$$s_{p}\Amalg s_{q}(z)=\begin{array}{ll}s_{p}(z)&{\rm{\ if\ }}\vert z\vert \leq 1,\\ (s_{q}(z^{-1}))^{-1}&{\rm{\ if\ }}\vert z\vert \geq 1.\end{array}$$
We are only considering, here, the case of $p=\frac{3}{7}$. It follows from Tan Lei's Theorem (\cite{TL, R1}) that $s_{3/7}\Amalg s_{q}$ is Thurston equivalent to a quadratic rational map (unique up to M\"obius conjugacy) if and only if $q{\rm{\ mod\ }}1\in (-\frac{1}{3},\frac{1}{3})$. As in the case of captures, it is essentially a folklore result  that, for $q_{1}\in (\frac{1}{7},\frac{2}{7})\cup (\frac{5}{7},\frac{6}{7})$, branched coverings $s_{3/7}\Amalg s_{q_{1}}$ and $s_{3/7}\Amalg s_{q_{2}}$ are Thurston equivalent if and only if $s_{q_{1}}=s_{q_{2}}$, that is, $e^{2\pi iq_{1}}$ and $e^{2\pi iq_{2}}$ are endpoints of the same minor leaf (so that $L_{q_{1}}=L_{q_{2}}$). It is likely that the corresponding result for $q_{1}\in (-\frac{1}{7},\frac{1}{7})$ is true, and can be proved from the analogous result for captures. It is hoped to include a proof of this in another paper in the near future. For $q_{1}\in (\frac{2}{7},\frac{1}{3})\cup (\frac{2}{3},\frac{5}{7})$, this is far from true, as the result of this paper shows.

There is a way of describing lamination maps $s_{q}$ and matings up to Thurston equivalence using homeomorphisms $\sigma _{\beta }$, which has been used in \cite {R1} and  \cite{R2} and which will be useful here, too. First we describe $s_{q}$ up to Thurston equivalence, for any odd denominator rational $q$ in $(0,1)$. Let $s_{0}$ denote the map $s_{0}(z)=z^{2}$. Suppose that $e^{2\pi i q}$ has period $r$ under $s_{0}$. Let $\beta $ be an arc from $0$ to  a point near $e^{2\pi i q}$ which does not intersect the unit circle. Let $\zeta $ be a similar path from $0$ to $s_{0}^{r-1}({\rm{end}}(\beta ))$, which is in a small neighbourhood of $s_{0}^{-1}(\beta )$. Then the Thurston equivalence class of 
$$\sigma _{\zeta }^{-1}\circ \sigma _{\beta }\circ s_{0}$$
is uniquely determined, assuming that sufficiently small neighbourhoods were chosen. We can choose $\beta $ to intersect the minor leaf of $L_{q}$ exactly once, and $\zeta $ so as not to intersect $L_{q}\cup S^{1}$ at all. It is then clear that $\sigma _{\zeta }^{-1}\circ \sigma _{\beta }\circ s_{0}$ preserves $L_{q}$ up to isotopy, and is Thurston equivalent to $s_{q}$. 

For the similar description of matings, we can describe $s_{p}\Amalg s_{q}$ as follows, again assuming that $e^{2\pi iq}$ is of period $r$ under $s_{0}$. We take $\beta $ to be a path from $\infty $ to a small neighbourhood of $e^{-2\pi iq}$ (not $e^{2\pi iq}$, because the identification with the unit disc is by $z\mapsto z^{-1}$) which does not cross the unit circle, and let $\zeta $ a path from $\infty $ to  a small neighbourhood of $s_{p}^{r-1}({\rm{end}}(\beta ))$ which is in a small neighbourhood of $s_{p}^{-1}(\beta )$. Then 
$$\sigma _{\zeta }^{-1}\circ \sigma _{\beta }\circ s_{p}$$
is Thurston equivalent to $s_{p}\Amalg s_{q}$.

\subsection{The Mating and Capture Theorems}\label{1.5}

We shall prove the following. 
\begin{captheorem} 
\begin{itemize}
\item There is a pair of equivalent captures by the aeroplane polynomial with the nonperiodic critical point of preperiod $13$, such that  the capture paths have distinct endpoints, but the $S^{1}$-crossing points of the capture paths are both in $\{ e^{2\pi it}:t\in (\frac{2}{7},\frac{9}{28})\} $.

\item For each integer $n\geq 0$ there is a set of  $n+2$ equivalent captures  by the aeroplane polynomial, with  the nonperiodic critical point of preperiod $26\cdot 2^{n}-11$, and such that all the capture paths have distinct endpoints, but the $S^{1}$-crossing points of the capture paths are all in $\{ e^{2\pi it}:t\in (\frac{2}{7},\frac{9}{28})\} $.\end{itemize}
 \end{captheorem}

\begin{mattheorem} 
For each integer $n\geq 0$ there is a set of  $n+2$ equivalent matings  $s_{3/7}\Amalg s_{q}$ with the aeroplane polynomial, with period $30.2^{n}-12$   for the second critical point,  and $q\in (\frac{19}{28},\frac{5}{7})$ for all such $q$. \end{mattheorem}

\subsection{Symbolic Dynamics}\label{1.4}

In order to identify points in $S^{1}$, and leaves and gaps of $L_{3/7}$, and, in order to do our calculations, we shall use the the symbolic dynamics for $s_{3/7}$ which is used in \cite{R5} and \cite{R6}, with letters $L_{j}$ and $R_{j}$ for $1\leq j\leq 3$, and $C$, $BC$ and $UC$.  The labels on regions are as shown. The vertical lines separating regions  join points $e^{\pm 2\pi i(r/14)}$ on the unit circle, for $1\leq r\leq 6$. We use $D(L_{j})$, $D(R_{j})$ and $D(C)$ to denote the closures of the regions labelled by $L_{j}$, $R_{j}$ or $C$. We define 
$$D(UC)=D(C)\cap \{ z:{\rm{Im}}(z)\geq 0\} ,$$
$$D(BC)=D(C)\cap \{ z:{\rm{Im}}(z)\leq 0\} .$$
\begin{figure}
\centering{\includegraphics[width=8cm]{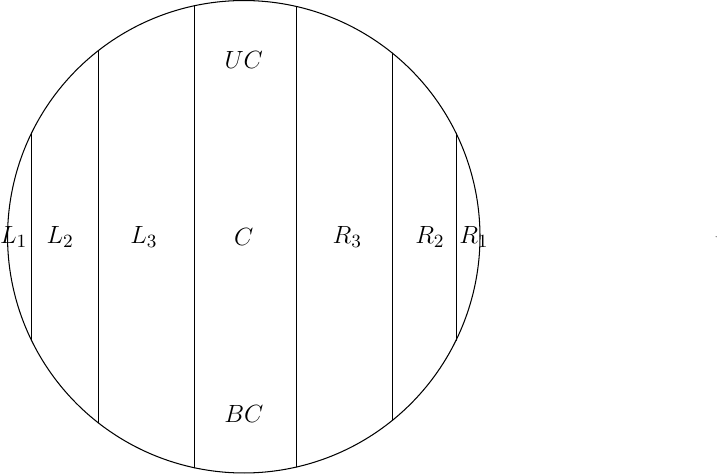}}
\caption{}
\end{figure}
Under the map $s_{3/7}$, regions map as follows:
$$\begin{array}{lll}D(R_{1}),\ D(L_{1})&\ \to\ &D(R_{1})\cup D(R_{2}),\\D(R_{2}),\ D(L_{2})&\to\ &D(R_{3})\cup D(C),\\ D(R_{3}),\ D(L_{3})\ &\to\ & D(L_{3})\cup D(L_{2}),\\ D(C)&\to\ & D(L_{1}).\end{array}$$
Each $D(R_{j})$ or $D(L_{j})$ is mapped homeomorphically by $s_{3/7}$ onto its image, while $D(C)$ is mapped with degree two onto $D(L_{1})$ by $s_{3/7}$. Each of $D(BC), D(UC)$ is also mapped onto $D(L_{1})$.

If $w=w_{0}\cdots w_{n}$ is a word in the letters $L_{j}$, $R_{j}$, $BC$, $UC$ and $C$ then we define
$$D(w)=\cap _{k=0}^{n}s_{3/7}^{-k}(D(w_{k})).$$
We can make a similar definition for infinite words, in which case $D(w)$ is a point on the unit circle $S^{1}$, or a leaf of $L_{3/7}$, or the closure of a gap of $L_{3/7}$. 
We define
$$Z_{m}=Z_{m}(s_{3/7})=s_{3/7}^{-m}(\{ 0,s_{3/7}(0),s_{3/7}^{2}(0)\} ),$$
and
$$Z_{\infty }=\cup _{m\geq 0}Z_{m}.$$
 If $z$ is the point of $Z_{\infty }$ of lowest preperiod in $D(w)$, then we write $z=p(w)$ and we sometimes write $D(z)$ for $D(w)$.

We shall often consider words which contain only the letters $L_{3}$, $L_{2}$ and $R_{3}$, possibly with $C$ at the end. For any such word $w$, $D(w)$ is bounded in the unit disc by vertical chords. For such words, we shall use the notation $v<w$ if $D(v)$ is to the right of $D(w)$ in the unit disc. We shall avoid this notation in the case when $v$ is a prefix of $w$ or vice versa. We shall also sometimes write $v\leq w$, meaning that $v<w$ or $v=w$. For such a word $w$, we shall write $\beta (w)$ for a capture  path ending at $p(w)$ and crossing into $D(w)$ on the upper unit circle. The exact crossing point does not matter, because all such paths are homotopic via homotopies which fix the forward orbit of $p(w)$.

\section{Proofs}\label{2}
\subsection{The basic idea}\label{2.1}
The following lemma is the basic idea. The rest of the proof of the capture and mating theorems is obtained by variations on this.
\begin{ulemma}
Define
$$u_{0}=L_{3}L_{2}R_{3}L_{3}L_{2}C,$$
$$v_{0}=L_{3}L_{2}R_{3}L_{3}(L_{2}R_{3})^{2},$$
$$w_{0}=L_{3}L_{2}R_{3}L_{3}^{5}.$$
 Then the captures $\sigma _{\beta (w_{0}u_{0})}\circ s_{3/7}$ and $\sigma _{\beta (v_{0}u_{0})}\circ s_{3/7}$ are Thurston equivalent.

\end{ulemma}

\noindent {\em{Proof.}} Write $s_{3/7}=s$ from now on. 
We omit the index $0$ and simply write $v$, $w$ and $u$ for $v_{0}$, $w_{0}$ and $u_{0}$. So take $x=p(wu)$ and $y=p(vu)$. Then $D(y)$ is to the right of $D(x)$, although both are in $D(L_{3}L_{2}R_{3}L_{3})$. Note, also, that $p(u)$  is in the forward orbit under $s$ of $x$ and $y$, and there is no other point in the forward orbit of $y$ which is between $D(y)$ and $D(x)$.   So now we take the subdisc $D'$ of the unit disc which is bounded by $D(x)$ and $D(y)$. It does not matter whether we choose $D'$ to contain $D(x)$ (or $D(y)$) or not, because we are only concerned with the isotopy class with respect to the strict forward orbit of $y$ but we will choose it to be strictly between $D(x)$ and $D(y)$. Choose a point $z\in Z_{\infty }$ whose word starts with $L_{3}$, and contains only the letters $L_{3}$, $L_{2}$ and $R_{3}$ such that $z$ is not in $\cup _{n=0}^{\infty }s^{-n}(D'\cup D(y)\cup D(x))$, and $D(z)$ is to the right of $D(y)$. We are now going to use $\sigma _{\beta (z)}\circ s$ as a base for the construction. The exact position of $z$ does not matter.

 Let $\alpha $ be the simple loop based at $z$, which crosses out of the gap containing $z$ on the upper unit circle, then crosses back to trace round $D'$ in an anti-clockwise direction, crosses out on the upper boundary of  the unit circle and then back across the upper unit circle into the gap containing $z$ and then ends at $z$. 
\begin{figure}
\centering{\includegraphics[width=8cm]{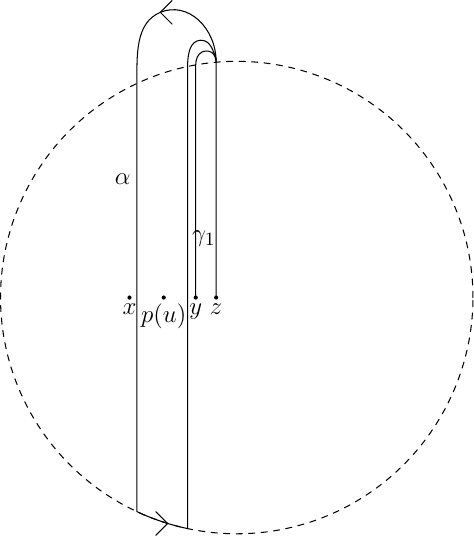}} 
\caption{$\alpha $}
\end{figure}
 Let 
 $$O(z)=\{ s^{m}(z):m\geq 0\} ,$$
and similarly for any other point replacing $z$. Note that $p(u)=s^{8}(y)=s^{8}(x)$ is of preperiod $5$ under $s$, and so $O(p(u))\subset (\sigma _{\beta (z)}\circ s)^{-5}(O(z))$. and for any $m\geq 5$, we have $s^{5-m}(O(p(u))\subset (\sigma _{\beta (z)}\circ s)^{-m}(O(z))$.  Now for any $m\geq 5$ we can find a homeomorphism 
  $$\psi_{m}:(\overline{\mathbb C},(\sigma _{\beta (z)}\circ s)^{-m}(O(z)))\to (\overline{\mathbb C},(\sigma _{\beta (z)}\circ s)^{-m}(O(z)))$$
  such that
   \begin{equation}\label{2.1.1} (\sigma  _{\beta (z)}\circ s,(\sigma _{\beta (z)}\circ s)^{-m}(O(z)))\simeq _{\psi _{m}}(\sigma _{\alpha }\circ \sigma _{\beta (z)}\circ  s,(\sigma _{\beta (z)}\circ s)^{-m}(O(z)))\end{equation}
where this notation $\simeq _{\psi _{m}}$ (which is also used in \cite{R1, R2, R3, R5, R6}) means that $\psi _{m}\circ \sigma _{\beta (z)}\circ s\circ \psi _{m}^{-1}$ and  $\sigma _{\alpha }\circ \sigma _{\beta (z)}\circ  s$ are homotopic via a homotopy through critically finite branched coverings which all take the same values on $(\sigma _{\beta (z)}\circ s)^{-m}(O(z))$. In particular, $\sigma _{\beta (z)}\circ s$ and $\sigma _{\alpha }\circ \sigma _{\beta (z)}\circ  s$ are Thurston equivalent. Of course, Thurston equivalence is an equivalence relation, while $\simeq _{\psi _{m}}$ is not. Nevertheless this notation gives useful information. Similarly to \cite{R5}, (3.3, 6.3 or 8.1) we can define $\psi _{m}$ by induction on $m$, setting  $\psi _{5}$ to be the identity. Then we define $\psi _{m+1}$ in terms of $\psi _{m}$ for $m\geq 5$ by 
$$\psi _{m}\circ \sigma _{\beta (z)}\circ s=\sigma _{\alpha }\circ \sigma _{\beta (z)}\circ  s\circ \psi _{m+1},$$
together with the property that $\psi _{m+1}$ is isotopic to $\psi _{m}$ via an isotopy which is the identity on $(\sigma _{\beta (z)}\circ s)^{-5}(O(z))$, which contains $O(p(u))$. It is then also true that $\psi _{m}$ and $\psi _{m+1}$ are isotopic via an isotopy which is the identity on $(\sigma _{\beta (z)}\circ s)^{-m}(O(z))$. Defining $\xi _{m}$ by
$$\psi _{m+1}=\xi _{m}\circ \psi _{m},$$
we have 
$$\sigma _{\beta (z)}\circ s=\sigma _{\alpha }\circ \sigma _{\beta (z)}\circ  s\circ \xi _{5},$$
and, for $m\geq 5$,
$$\sigma _{\alpha }\circ \sigma _{\beta (z)}\circ s\circ \xi _{m+1}=\xi _{m}\circ \sigma _{\alpha }\circ \sigma _{\beta (z)}\circ s.$$
We then have
$$\psi _{m}=\xi _{m-1}\circ \cdots \circ \xi _{5}.$$
We need $\psi _{m}$ for $m$ equal to the preperiod of $y$ (so $m=13 $ in this case, which is also the preperiod of  $x$), up to isotopy preserving $O(y)\subset s^{5-m}(O(p(u)))\subset (\sigma _{\beta (z)}\circ s)^{-m}(O(z))$. For such an $m$, we write $\psi _{m}=\psi $. Now each $\xi _{k}$ is a composition of disc- exchanging involutions. We see this first for $\xi _{5}=\psi _{6}$, which is a single involution  whose support is an arbitrarily small neighbourhood of
$$s^{-1}(\beta (z)\cup \alpha \cup D'),$$
that is, approximately the union of $s^{-1}(D')$ and the diagonal in the exterior of the unit disc connecting these. The involution  interchanges the components of $s^{-1}(D')$ --- tracing a path halfway round the boundary of the supporting disc  in a {\em{clockwise}} direction. This ensures that $\xi _{5}^{-1}=\psi _{6}^{-1}$ exchanges discs, tracing a path halfway round the boundary of the supporting disc in an {\em{anticlockwise}} direction. Similarly we have
$${\rm{supp}}(\xi _{k+1})=(\sigma _{\alpha }\circ \sigma _{\beta (z)}\circ s)^{-1}({\rm{supp}}(\xi _{k})).$$
We only need to consider the components of the support of $\xi _{k}$ which intersect $O(y)$. Note that the support of $\xi _{5}$ does not intersect $\alpha $, although it does intersect $\beta (z)$ near $\infty $. We shall see inductively that, for all $6\leq k\leq 12$, the components of $\xi _{k}$ which intersect $O(y)$ do not intersect $\alpha \cup \beta (z)$, but do intersect $O(x)$. So we have
$${\rm{supp}}(\xi _{6})=(\sigma _{\beta (z)}\circ s)^{-1}({\rm{supp}}(\xi _{5})),$$
and for $6\leq k\leq 12$,
$${\rm{supp}}(\xi _{k+1})=s^{-1}({\rm{supp}}(\xi _{k})).$$
   We can describe the involutions using the symbolic dynamics. The involution of $\xi _{5}$ is described by
$$L_{3}D'{\rm{(bot)}}\leftrightarrow {\rm{(top)}}R_{3}D'.$$
\begin{figure}
\centering{\includegraphics[width=8cm]{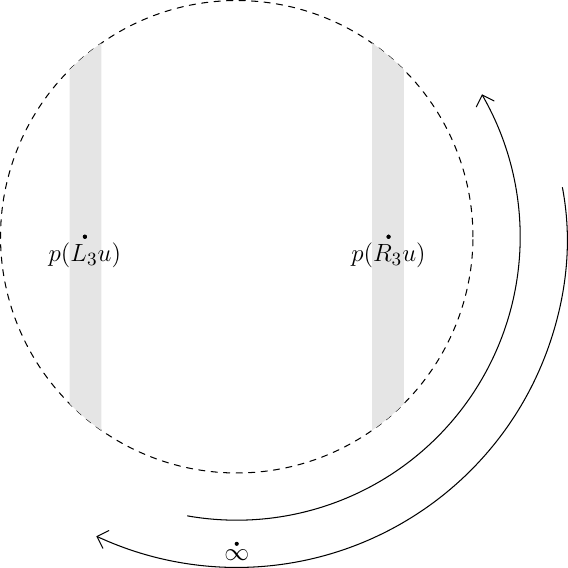}}
\caption{$\xi _{5}$}
\end{figure}
The support of $\xi _{5}$ thus contains just the points $p(L_{3}u)$ and $p(R_{3}u)$ in $O(x)$ and $O(y)$ respectively.  We have the following components in the support of $\xi _{k}$ for $5\leq k\leq 12$. In each case, two discs $u_{1}D'$ and $u_{2}D'$ are connected by an arc in the exterior of the unit disc. The labels ``top'' and ``bot'' indicate whether the arc joins onto the disc $u_{j}D'$ on the upper or lower half of the unit circle.
We then only need to consider the components of the support of $\xi _{k}$ which intersect $O(y)$. Now $p(u)$ is the only point in $D'$ in the forward orbit of $y$  (and $x$), and it turns out that a single component of the support of $\xi _{k}$ contains the points in the backward orbit of both $x$ and $y$ which map forward under $s^{k+1}$ to $u$. For $6\leq k\leq 11$, the action of $\xi _{k}$ is represented by: 
$$L_{3}^{2}D'{\rm{(top)}}\leftrightarrow{\rm{(bot)}}L_{2}R_{3}D',$$
$$L_{3}^{3}D'{\rm{(bot)}}\leftrightarrow{\rm{(bot)}}R_{3}L_{2}R_{3}D',$$
$$L_{3}^{4}D'{\rm{(top)}}\leftrightarrow{\rm{(top)}}L_{2}R_{3}L_{2}R_{3}D',$$
$$L_{3}^{5}D'{\rm{(bot)}}\leftrightarrow{\rm{(bot)}}L_{3}L_{2}R_{3}L_{2}R_{3}D',$$
$$R_{3}L_{3}^{5}D'{\rm{(bot)}}\leftrightarrow{\rm{(bot)}}R_{3}L_{3}L_{2}R_{3}L_{2}R_{3}D',$$
$$L_{2}R_{3}L_{3}^{5}D'{\rm{(top)}}\leftrightarrow{\rm{(top)}}L_{2}R_{3}L_{3}L_{2}R_{3}L_{2}R_{3}D',$$
and finally, for $k=12$:
$$vD'{\rm{(bot)}}\leftrightarrow {\rm{(bot)}}wD'.$$ 

The support of this last component contains both $x$ and $y$, and interchanges them. Define the path $\gamma _{1}$ by
$$\beta (y)=\beta (z)*\gamma _{1},$$ 
where $*$ denotes the usual composition of paths. So the endpoints of $\gamma _{1}$ are $z$ and $y$. We may take this as a pointwise definition of $\gamma _{1}$, but we are only interested in $\gamma _{1}$ up to homotopy preserving $O(y)$ and $O(z)$. Up to such a homotopy, the path $\gamma _{1}$ is an arc which crosses out of the gap containing $z$ across the upper unit circle then into the gap containing $y$, again across the upper unit circle, and ends at $y$. Then we write 
$$\gamma _{2}=\alpha *\psi (\gamma _{1}).$$
 Then
$$\sigma _{\beta (y)}\circ s\simeq _{\rm{identity}}\sigma _{\gamma _{1}}\circ \sigma _{\beta (z)}\circ s\simeq _{\psi }\sigma _{\gamma _{2}}\circ \sigma _{\beta (z)}\circ s.$$

 None of the involutions in the support of $\psi $ intersects $\gamma _{1}$, except for the one in the support of $\xi _{12}$  which interchanges the discs $D(x)$ and $D(y)$, with connecting arc along the lower unit circle. So the endpoints of $\gamma _{2}$ are $z$ and $x$. Because this disc exchange (like all of them) is in a clockwise direction,  $\psi (\gamma _{1})$ passes to the left of $D(y)$, while $\alpha $ passes to the right of $vD_{1}'$ where $D_{1}'\subset D'$ is defined by $vD_{1}'=vD'\setminus D'$. 
 However $vD_{1}'\cap O(x)=\emptyset $.  Note that when discs are exchanged, upper boundaries go to upper boundaries, and lower boundaries to lower. So $\psi (\gamma _{1})$ hooks round $wD_{1}'$, but  $wD_{1}'\cap O(x)=\emptyset $ --- because $vD_{1}'\cap O(y)=\emptyset $. So 
 $$\beta (z)*\gamma _{2}=\beta (z)*\alpha *\psi (\gamma _{1})=\beta (x){\ \rm{rel\ }}O(x),$$
 and 
\begin{equation}\label{2.1.2}\sigma _{\gamma _{2}}\circ \sigma _{\beta (z)}\circ s\simeq _{\psi }\sigma _{\beta (x)}.\circ s\end{equation}
\Box

\subsection{Extensions of \ref{2.1}}\label{2.2}

Write
$$a=L_{3}(L_{2}R_{3})^{2},\ \ b=L_{3}^{5}.$$
It is clear that the pair $(a,b)$ plays a key role in \ref{2.1}. Note that
$$v_{0}=L_{3}L_{2}R_{3}a,\ \ w_{0}=L_{3}L_{2}R_{3}b.$$
The following properties hold.
\begin{itemize}
\item[1.] The words $v_{0}$ and $w_{0}$ contain only the letters $L_{3}$, $L_{2}$ and $R_{3}$, and $u_{0}$ contains only the letters $L_{3}$, $L_{2}$ and $R_{3}$ apart from the last letter $C$.
\item[2.]  We have 
$$v_{0}<u_{0}<w_{0}$$
and $u_{0}$ is the only suffix $u$ of $v_{0}u_{0}$ which satisfies 
$$v_{0}u_{0}<u<w_{0}u_{0},$$
and $w_{0}u_{0}$ is obtained from $v_{0}u_{0}$ by changing the word $a$ preceding $u_{0}$ in $v_{0}u_{0}$ to $b$.
\end{itemize}

These are essentially the only properties of $u_{0}$, $v_{0}$ and $w_{0}$ which  are used in the proof of \ref{2.1}. We want words $v_{0}$ and $w_{0}$ such that $D(v_{0})$ and $D(w_{0})$ are in the subset of the unit disc bounded by the leaves of $L_{3/7}$ with endpoints at $e^{\pm 2\pi i(2/7)}$ and $e^{\pm 2\pi i(1/3)}$. This will be true provided that the common prefix starts  with $L_{3}^{2i+1}L_{2}$ for some $i\geq 0$. 

This motivates the following. We choose the words $c$ and $d$ to minimise the preperiods and periods in the Capture and Mating Theorems for all $n\geq 1$. Of course other choices are possible.
Let $a$ and $b$ be as above and let $v_{0}$, $w_{0}$ and $u_{0}$ be as in \ref{2.1}. In addition define
$$c=L_{3}^{3}L_{2}C,$$
$$d=L_{3}^{2}a,$$
$$u_{0}'=L_{3}L_{2}R_{3}c,$$
$$t_{0}=L_{3}L_{2}R_{3}d,$$
and define $v_{k}$, $w_{k}$, $t_{k}$ and $u_{k}$ inductively for $k\geq 1$ by
$$v_{k}=v_{k-1}t_{k-1}a,$$
$$w_{k}=v_{k-1}t_{k-1}b,$$
$$u_{k}=v_{k-1}t_{k-1}c,$$
$$t_{k}=v_{k-1}t_{k-1}d.$$

We write $\vert v\vert $ for the number of letters in a word $v$. Note that, for $k\geq 1$,
$$\vert v_{k}\vert +\vert t_{k}\vert =2(\vert v_{k-1}\vert +\vert t_{k-1}\vert )+12,$$
$$\vert v_{k}\vert -\vert t_{k}\vert =2(\vert v_{k-1}\vert -\vert t_{k-1}\vert )-2,$$
and so
$$\vert v_{n}\vert =2\vert v_{n-1}\vert +5=13.2^{n}-5,$$
$$\vert t_{n}\vert =2\vert t_{n-1}\vert +7=17.2^{n}-7,$$
and for $n\geq 1$ (not $n=0$)
$$\vert u_{n}\vert =13.2^{n}-6.$$
So we have
$$\vert v_{n}\vert +\vert u_{n}\vert =26.2^{n}-11$$
for all $n\geq 1,$ and 
$$\vert v_{n}\vert +\vert t_{n}\vert =30.2^{n}-12$$
for all $n\geq 0$. This will give the numbers in the Capture and Mating Theorems.

\begin{lemma}\label{2.3} 
 $v_{k}$, $u_{k}$ and $w_{k}$ satisfy the conditions 1---2 of \ref{2.2}.\end{lemma}

\noindent {\em{Proof}}

Each of $a$, $b$, $c$, $d$,  and $u_{0}$ has an odd number of letters $L_{3}$ and $L_{2}$.  Hence, the same is true for $v_{0}$, $w_{0}$ and $t_{0}$, which are obtained by prefixing $a$, $b$ and $d$ by $L_{3}L_{2}R_{3}$. Hence, by induction,  the same is true of $v_{k}$, $w_{k}$, $u_{k}$ and $t_{k}$ for all $k\geq 0$. Since the words $v_{k}$, $w_{k}$, $u_{k}$ and $t_{k}$ for $k>0$ are obtained from $a$, $b$, $c$ and $d$ by adjoining the common prefix $v_{k-1}t_{k-1}$, which has an even number of letters $L_{3}$ and $L_{2}$, and since 
$$a<d<c<b,$$
 we have
 $$v_{0}<t_{0}<u_{0}'<w_{0},$$
 and,  for all $k>0$,
\begin{equation}\label{2.3.1}v_{k}<t_{k}<u_{k}<w_{k}.\end{equation}
 Then, for $k>1$, since
$t_{k-1}<u_{k-1}$ we have $v_{k-1}u_{k-1}<v_{k-1}t_{k-1}$,  
and
$$v_{k-1}u_{k-1}<v_{k}u_{k}.$$
Now the only condition of \ref{2.2} that needs to be checked is that $u_{k}$ is the only suffix  $u$ of $v_{k}u_{k}$ satisfying $v_{k}u_{k}<u<w_{k}u_{k}$.  Since $v_{k-1}t_{k-1}$ is a common prefix of $v_{k}$ and $w_{k}$, we only need to consider proper suffixes of $v_{k}u_{k}$ which start with $v_{k-1}t_{k-1}$. To do this, we prove something more general about occurrences of $v_{k-1}t_{k-1}$.

\begin{lemma}\label{2.4}The only  occurrences of $v_{k-1}t_{k-1}$ in
\begin{equation}\label{2.4.1}v_{k}t_{k}, t_{k}av_{k},\ t_{k}au_{k}\end{equation}
 are those visible in the definitions of $v_{k}$, $t_{k}$ and $u_{k}$ above, that is, at the start of these words. \end{lemma}
 
\noindent {\em Proof} We prove this by induction on $k$. We can easily check that it is true for $k=1$, by inspection. Now suppose that it is true for $k-1$. Now $v_{k-1}t_{k-1}$ starts with $v_{k-1}$, which starts with $v_{k-2}t_{k-2}$. Now
 $$\begin{array}{l}
 v_{k}t_{k}=v_{k-1}t_{k-1}av_{k-1}t_{k-1}d,\\ 
 t_{k}av_{k}=v_{k-1}t_{k-1}dat_{k-1}v_{k-1}a,\\
  t_{k}au_{k}=v_{k-1}t_{k-1}dav_{k-1}t_{k-1}c.\end{array}$$
  So any occurrence of $v_{k-1}t_{k-1}$ in any of the words of (\ref{2.4.1}) contains a $v_{k-2}t_{k-2}$ which must be in one of the corresponding words for $k-2$ and must coincide with the start of one of the words $v_{k-1}$, $t_{k-1}$ or $u_{k-1}$. Then by inspection we see that the only possibility for $v_{k-1}t_{k-1}$  in any of the words of (\ref{2.4.1}) is at the start of one of the  subwords $v_{k}$, $t_{k}$ or $u_{k}$.

\Box

We then have the following corollary.

\begin{theorem}\label{2.5} For $v_{k}$, $w_{k}$ and  $u_{k}$ defined as above, the captures $\sigma _{\beta (w_{k}u_{k})}\circ s_{3/7}$ and $\sigma _{\beta (v_{k}u_{k})}\circ s_{3/7}$ are Thurston equivalent.
\end{theorem}

\noindent {\em{Proof}} Define $y=p(v_{k}u_{k})$ and $x=p(w_{k}u_{k})$ similarly to \ref{2.1}, and define $\alpha $, $D'$, $\gamma _{1}$ and $\gamma _{2}$ as in \ref{2.1}. As before, choose $z$ to the right of $D(y)$ and to the left of $D(C)$, and outside $\cup _{n\geq 0}s^{-n}(D'\cup D(x)\cup D(y))$. Let $m=\vert v_{k}\vert =13\cdot 2^{k}-5$. Let $S$ be the local inverse of $s^{m}$ defined by $v_{k}$, so that $Sp(u_{k}))=x$. Then $S(D'\cup D(x)\cup D(y))$ intersects  $D'\cup D(x)\cup D(y)$, but, by \ref{2.3}, it does not contain any points of $O(x)$ apart from $x$ itself. Then we can choose $z$ arbitrarily close to $\cup _{n\geq 0}S^{n}(D'\cup D(x)\cup D(y))$, and on the right of this set. Then define the sequences of homeomorphisms $\psi _{m}$ and $\xi _{\ell }$ as in \ref{2.1}. We need to consider the effect of disc exchanges in the support of $\xi _{\ell }$ on $\gamma _{1}$, as in \ref{2.1}. We only need to consider those which pull back to ones which  intersect $\gamma _{1}$, and which also intersect $O(y)$. The first few disc exchanges are as in \ref{2.1}. The rest are all of the form $u'aD'\leftrightarrow u'bD'$, where $u'$ is a prefix of $u_{k}$, and, also, $u'au_{k}$ is a suffix of $v_{k}u_{k}$, that is, $u'a$ is a suffix of $v_{k}$. By the conditions imposed on $z$, it is never in either of the sets $u'aD'$ or $u'bD'$, and is only strictly between them if $D'$ is also. So $\gamma _{1}$ does not intersect the disc exchange given by $u'aD'\leftrightarrow u'bD'$, if $y$ is not in $u'aD'\cup u'bD'$, and $y$ is  not between $u'aD'$ and $u'bD'$. So there is no intersection, except possibly when $u'a$ is a suffix of $v_{k}=v_{k-1}t_{k-1}a$,  and $u'$ is a prefix of $v_{k}$. So now we consider when this can happen. Take  the largest $i$ such that $v_{i}t_{i}$ is a (not necessarily proper) prefix of $u'$, and hence also a prefix of $v_{k}$. Then by \ref{2.4},  the only occurrences of $v_{i}t_{i}$ in $v_{k-1}t_{k-1}$ are those in the subwords $v_{i+1}$ and $t_{i+1}$ which occur in the construction of $v_{k-1}$ and $t_{k-1}$. So the only possibility is that $u'=v_{k-1}t_{k-1}$ and $u'a=v_{k}$. This gives the disc exchange $v_{k}D'\leftrightarrow w_{k}D'$, which has exactly the right effect on $\gamma _{1}$, as in \ref{2.1}.

\Box

To obtain equivalences between more than two captures, we need a slight generalisation. 
For $n>k$, define $w_{k,n}$ and $u_{k,n}$ and $t_{k,n}$ by replacing every occurrence of $a$ which precedes a word $v_{k}$ in $v_{n}$ -- or $u_{n}$ -- by $b$. Note that by \ref{2.4},  applied to $v_{k-1}t_{k-1}$ every occurrence of  $v_{k}$ in these longer words is indeed preceded by $a$, since $d=L_{3}^{2}a$.

\begin{lemma}\label{2.6} The suffixes $u$ of $v_{n}u_{n}$ which satisfy $v_{n}u_{n}<u<w_{k,n}u_{k,n}$ are precisely those which start with $v_{k}$ or $t_{k}$.\end{lemma}

\noindent{\em{Proof}} We consider the case $k=0$ since the proofs are exactly similar. It is clear that we only need to consider prefixes which start with $v_{0}$ or $t_{0}$. Every suffix $u$ which starts with $t_{0}$ certainly satisfies $v_{n}<u<w_{0,n}$, because $w_{0,n}$ starts with $w_{0}$ and $v_{n}$ starts with $v_{0}$. So there is only a question about suffixes which start with $v_{0}$, and all of these start with $v_{0}t_{0}$, and hence with either $v_{1}=v_{0}t_{0}a$ or with $t_{1}=v_{0}t_{0}d$ --- or with $u_{1}=v_{0}t_{1}c$ if $n=1$. So then we are reduced to considering suffixes which start with $v_{1}$. Continuing in the same way, we see that all suffices which start with $v_{r}$ for any $1\leq r<n$ or with $u_{n}$ have the required property.\Box

The following completes the proof of the theorem for captures.

\begin{theorem}\label{2.7}  $\sigma _{\beta (w_{k,n}u_{k,n})}\circ s_{3/7}$ and $\sigma _{\beta (v_{n}u_{n})}\circ s_{3/7}$ are Thurston equivalent. \end{theorem}

\noindent{\em{Proof}} We consider $n>k\geq 0$ since \ref{2.1} covers the case $k=n=0$, and \ref{2.5} covers $k=n$ in general.
We consider the case $k=0$ since the other cases are exactly similar for any $k<n$.  By \ref{2.6}, the points in the forward orbit of $p(v_{n}u_{n})$ that are between $p(v_{n}u_{n})$ and $p(w_{0,n}u_{n})$ are precisely  the points $p(u)$ for all suffixes $u$ of $v_{n}u_{n}$ which start with $v_{0}$ or $t_{0}$.

Once again, we proceed much as in \ref{2.1} --- and \ref{2.5}. Similarly to there, we define
$$x=p(w_{0,n}u_{0,n}),\ \ y=p(v_{n}u_{n}).$$
 We define $D'$ to be the region strictly between $D(w_{0,n})$ and $D(v_{n}u_{n})$. Choose $z$ as in \ref{2.5}, but with $S$ being the local inverse defined by $v_{k}$. This is still possible, using \ref{2.6} rather than \ref{2.3}. As in \ref{2.1}, we let $\alpha $ be the  anticlockwise loop based at $z$ round $D'$, we define $\gamma _{1}$ by $\beta (z)*\gamma _{1}=\beta (y)$, and proceed to define $\xi _{k}$ and $\psi _{m}$ as in \ref{2.1}. As in \ref{2.1}, we want to prove that, for $\psi =\psi _{m}$ and $m=\vert v_{n}\vert +\vert u_{n}\vert -1=26.u^{n}-11$, 
 \begin{equation}\label{2.7.1}\beta (x)=\beta (z)*\alpha *\psi (\gamma _{1}){\rm{\ rel\ }}O(x).\end{equation}
 
 By the same argument as in \ref{2.5},  the only disc exchanges $u'aD'\leftrightarrow u'bD'$ which have to be considered are those for which $O(y)\cap u'aD'\neq \emptyset $ and $y\in D(u')$. By \ref{2.6} these are those for which $u'aw$ is a suffix of $v_{n}u_{n}$, where $w$ is a suffix of $v_{n}u_{n}$ which starts with $v_{0}$ or $t_{0}$, and $u'$ is also a prefix of $v_{n}u_{n}$. Using \ref{2.4}, we see that the possibilities are $u'aw$ must be a suffix of $v_{n}u_{n}$ which starts with $v_{i}$ or $t_{i}$ for some $i>0$, or (if $n=1$) with $u_{1}$. All such disc exchanges are properly contained in $D'$, and do not intersect $\alpha \cup \beta (z)$. 
 
 This time, in contrast to \ref{2.1} and \ref{2.5}, there are many disc exchanges intersecting $D'$  and $\gamma _{1}$,  along the lower unit circle and the upper unit circle.  The occurrences of $v_{0}t_{0}$ in $v_{n}u_{n}$ are separated by strings of the form $d^{r}a$ for a variable integer $r\geq 0$. The number of letters $L_{3}$ and $L_{2}$ in $d^{r}a$ is odd or even, depending on whether $r$ is even or odd, but each of $v_{0}$ and $d_{0}$ has an odd number of letters $L_{3}$ and $L_{2}$, and, apart from the very first occurrence of $v_{0}$, each occurrence of $v_{0}$ or $t_{0}$ in $v_{n}u_{n}$ is always preceded by $a$.  Now every exchange $\xi _{j}$ along the lower unit circle introduces a double hook in $\psi _{j+1}(\gamma _{1})$, and  every exchange $\xi _{j}$ along the upper unit circle is such that $\psi _{j+1}(\gamma _{1})$  and $\psi _{j}(\gamma _{1})$ have the same number of hooks, and the last hook on $\psi _{j+1}(\gamma _{1})$ is a reduction of the last hook  on $\psi _{j}(\gamma _{1})$. The picture shows $\psi _{\ell _{1}}(\gamma _{1})$, where $\ell _{1}$ is the first $\ell $ for which $\psi _{\ell }(\beta (y))\neq \beta (y)$, so that $\ell _{1}=\vert v_{k}\vert =13.2^{k}-5$ (see \ref{2.2}). Note that  $\psi _{\ell _{1}}(y)$ is to the right of $x$. If $\ell _{2}$ is the first $\ell >\ell _{1}$ for which $\psi _{\ell }(\gamma _{1})\neq \psi _{\ell _{1}}(\gamma _{1})$, so that $\ell _{2}=13.2^{k+1}-5$, then $\psi _{\ell _{2}}(y)$ is to the right of $x$, and the last hook is round $x$.
 \begin{figure}
\centering{\includegraphics[width=12cm]{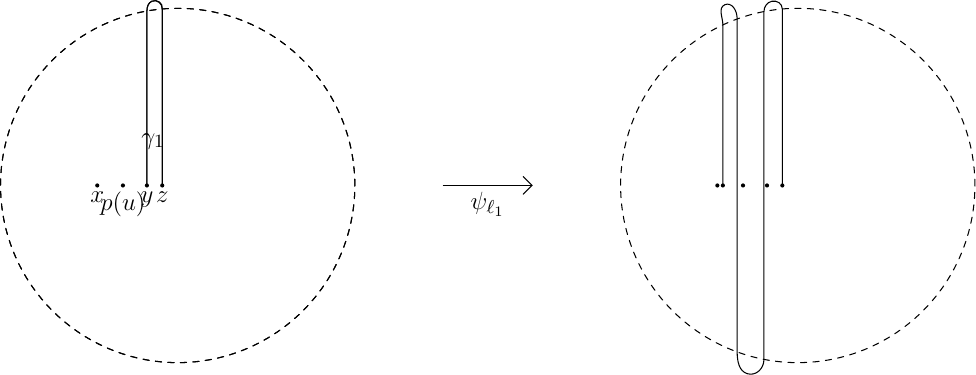}}
\caption{$\psi _{\ell _{1}}(\gamma _{1})$}
\end{figure}
In order to obtain (\ref{2.7.1}), we need the upward-pointing hooks on $\psi (\gamma _{1})$, that is, those which have the hook bridging the upper unit circle, to be homotopically trivial relative to $O(x)$.  This is true, because all the upward pointing hooks are  round subsets of $\psi (D_{1}')$ where (as in \ref{2.1}) $D_{1}'=v_{0}D'\setminus D'$ is to the right of $D'$, and  $D_{1}'\cap O(y)=\emptyset $ by \ref{2.6}. Therefore $\psi (D_{1}')\cap O(x)=\emptyset $. So we do have (\ref{2.7.1}), and, as in \ref{2.1}, this gives the desired equivalence.\Box
 
\subsection{Proof of the Matings Theorem}\label{2.8}
 
The proof of the Matings Theorem is then completed by the following lemma.
 
\begin{ulemma}
We continue with the same notation. In addition let $r_{k,n}$ be obtained from $t_{n}$ by replacing every occurrence of $a$ in $t_{n}$ before an occurrence of $v_{k}$ or $t_{k}$ in  $t_{n}$ by $b$.
Let $e^{2\pi iq_{k,n}}$ be the lower end of the leaf of $L_{3/7}$ coded by
$$(v_{n}t_{n})^{\infty }$$
if $k=n+1$, and by
$$(w_{k,n}r_{k,n})^{\infty }$$
 if $0\leq k\leq n$. Then all the matings  $s_{3/7}\Amalg s_{q_{k,n}}$ are Thurston equivalent. 
\end{ulemma}

\noindent{\em{Proof.}} 
 We start by considering the case $k=0=n$. Write $q_{0,0}=q_{0}$ and $q_{1,0}=q_{1}$ Note that $r_{0}=t_{0}$. As in \ref{2.1}, we drop the index $0$ and write  
 $$v=v_{0}=L_{3}L_{2}R_{3}L_{3}(L_{2}R_{3})^{2},$$
 $$u=u_{0}'=L_{3}L_{2}R_{3}L_{3}^{3}L_{2}C,$$ 
 $$w=w_{0}=L_{3}L_{2}R_{3}L_{3}^{5},$$ 
 $$t=t_{0}L_{3}L_{2}R_{3}L_{3}^{3}(L_{2}R_{3})^{2}.$$
   Note that the common length of $vt$ and $wt$ is $18$, and so in this particular case we are establishing an equivalence between matings with the aeroplane for which the second critical point has period $18$. 

We define 
$$x=p(wtwu),$$
and 
$$y=p(vtvu).$$
Then $x$ and $y$ have preperiod $33$.  As before, we write $s_{3/7}=s$. 

Choose a path $\zeta (y)$ from $\infty $ to $p(R_{3}vu)$ which crosses the unit circle just once, in the upper boundary of the gap of $L_{3/7}$ containing $p(R_{3}vu)$, then stays in this gap, and ends at $p(R_{3}vu)$. Then, using the description of 
$s\Amalg s_{q_{1}}$ up to Thurston equivalence given at the end of \ref{1.3},
$$s\Amalg s_{q_{1}}\simeq \sigma _{\zeta (y)}^{-1}\circ \sigma _{\beta (y)}\circ s.$$
This is because the gap of $L_{3/7}$ containing the  point $s^{j}(y)$, for $0\leq j\leq 17$, is closer to $\exp (2\pi i 2^{j}q_{1})$ than any other point $s^{i}(y)$ for $0\leq i\leq 17$ . Defining $\gamma _{1}$ by $\beta (z)*\gamma _{1}=\beta (y)$ as before, we also have 
$$ \sigma _{\zeta (y)}^{-1}\circ \sigma _{\beta (y)}\circ s\simeq  _{{\rm{identity}}}\sigma _{\zeta (y)}^{-1}\circ \sigma _{\gamma _{1}}\circ \sigma _{\beta (z)}\circ s.$$
 Similarly, define $\zeta (x)$ ending at $p(R_{3}wu)$, so that we have 
$$s\Amalg s_{q_{0}}\simeq \sigma _{\zeta (x)}^{-1}\circ \sigma _{\beta (x)}\circ s.$$

We are going to proceed much as in previous results, in particular, as in \ref{2.1}. So let $D'$ be the subset of the unit disc between $D(x)$ and $D(y)$, and let $z\in Z_{\infty}$ be such that $D(z)$ is to the right of $D(x)$,  and 
$$O(z)\cap (D'\cup D(x)\cup D(y))=\emptyset .$$
 Let $\alpha $ be the anticlockwise loop based at $z$ and  round $D'$. As in \ref{2.1}, we have a homeomorphism $\psi =\psi _{33}$  such that 
$$(\sigma _{\beta (z)}\circ s,(\sigma _{\beta (z)}\circ s)^{-33}(O(z)))\simeq _{\psi }(\sigma _{\alpha }\circ \sigma _{\beta (z)}\circ s,(\sigma _{\beta (z)}\circ s)^{-33}(O(z))).$$
 We have 
$$\sigma _{\zeta (y)}^{-1}\circ \sigma _{\gamma _{1}}\circ \sigma _{\beta (z)}\circ s \simeq _{\psi }\sigma _{\psi (\zeta(y) )}^{-1}\circ \sigma _{\psi (\gamma _{1})}\circ \sigma _{\alpha }\circ \sigma _{\beta (z)}\circ s.$$
Now the postcritical set $P(y)$ of  $\sigma _{\zeta (y)}^{-1}\circ \sigma _{\beta (y)}\circ s$ is 
$$\{ 0,s(0),s^{2}(0)\} \cup \{ s^{i}(y):0\leq y<18\} .$$
Note that it does not include $p(vu)$, and similarly the postcritical set $P(x)$ of  $\sigma _{\zeta (x)}^{-1}\circ \sigma _{\beta (x)}\circ s$ does not include $p(wu)$. So to prove the lemma, we only need to prove that
\begin{equation}\label{2.8.1}\beta (z)*\alpha *\psi (\gamma _{1})=\beta (x){\rm{\ rel\ }}P(x),\end{equation}
\begin{equation}\label{2.8.2} \psi (\zeta (y))=\zeta (x){\rm{\ rel\ }}P(x).\end{equation}
The proof of (\ref{2.8.1}) is much as in \ref{2.1} and \ref{2.7}. The support of $\xi _{\ell }$ intersects $\psi _{\ell }(\beta (y))$ only for $\ell =14$ and $\ell =32$. In both cases, just one component of the support of $\xi _{\ell }$ intersects $\psi _{\ell }(\beta (y))$, and restricted to this component, $\xi _{\ell }$ exchanges two discs. The homeomorphism $\xi _{14}$ interchanges a disc containing $p(vu)$ and $p(vtvu)$ with a disc containing $p(wu)$ and $p(wtvu)$, with $p(vu)$ and $p(vtvu)$ mapped to $p(wu)$ and $p(wtvu)$ respectively, and $\xi _{32}$ interchanges a disc containing $p(wtvu)$ with a disc containing $p(wtwu)$.   Figure 6 shows $\alpha *\psi _{15}(\beta (y))$ and $\alpha *\psi _{33}(\beta (y))$. The endpoint of $\psi _{15}(\beta (y))$ is $p(wtvu)$ and the endpoint of $\psi _{33}(\beta (y))$ is $p(wtwu)=x$.  Both of these paths loop round $p(wu)$ as shown (but because of space constraints, this point is not labelled). Since $p(wu)=s^{18}(p(x))$ is not in $P(x)$, this  does not matter, and we have  (\ref{2.8.1}) in the case $n=0$.
\begin{figure}
\centering{\includegraphics[width=12cm]{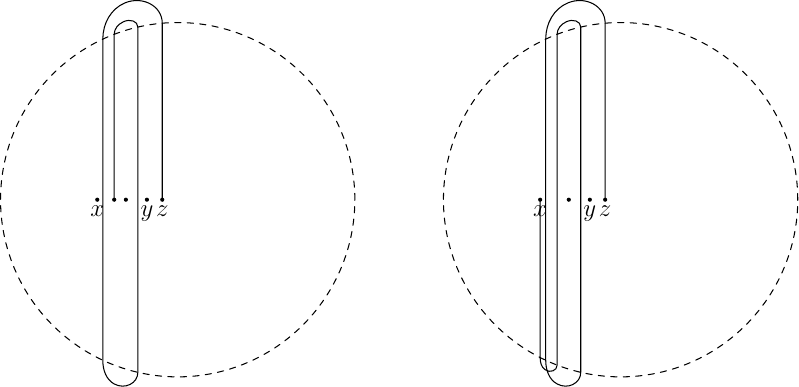}}
\caption{$\alpha *\psi _{15}(\beta (y))$ and $\alpha *\psi _{33}(\beta (y))$}
\end{figure}

Now we consider (\ref{2.8.2}).
 The action of $\xi _{7}=\psi _{8}$ on $\zeta (y)$ is as shown in Figure 7. Note that the support of $\xi _{7}$ intersects $\zeta (y)$ but does not contain $p(R_{3}vu)$ because $p(vu)\not\in D'$, because we chose $u$ with $t<u$, and therefore $vu<vt$.
\begin{figure}
\centering{\includegraphics[width=12cm]{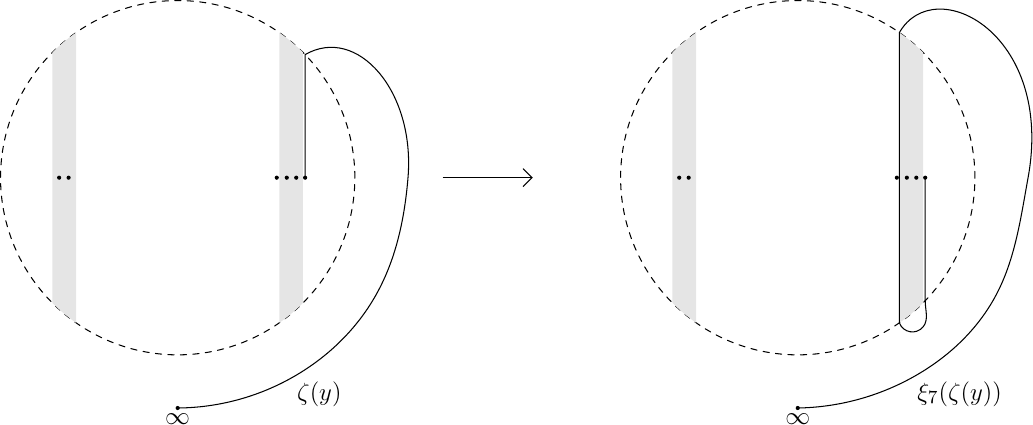}}
\caption{}
\end{figure}
The two discs in the support of $\xi _{7}$ interchange in a clockwise direction round $\infty $, and the action on $\zeta (y)$ is such that $\xi _{7}(\zeta (y))$ hooks round the nearby component of $s^{-1}(D')$. The points $p(L_{3}u)$ and $p(R_{3}u)$ in $s^{-1}(D')$ are interchanged by $\xi _{7}$. Then there is a composition of disc exchanges similar to those in \ref{2.1}. The support of $\xi _{\ell}$ does not intersect $\psi _{8}(\zeta (y))$ for $7<\ell <14$ but does intersect it for $\ell =15$. The endpoint of $\psi _{16}(\zeta (y))$ is then at $p(R_{3}wu)$. For $\ell \geq 16$, the discs in the support of $\xi _{\ell }$ contain only points of preperiod higher than $16$, and $16 $ is the preperiod of $p(R_{3}wu)$. The discs are sufficiently small that any such disc which intersects $\psi _{16}(\zeta (y))$ does not intersect $P(x)$. In particular, $\psi _{33}(\zeta (y))$ is homotopic to $\psi _{16}(y)$ relative to $P(x)$. This gives (\ref{2.8.2}) in the case $k=0$. The case $k=n$, for any $n\geq 0$, of (\ref{2.8.1}) and (\ref{2.8.2}), is exactly similar.

So now we consider the case $n>k\geq 0$. Write 
$$y=p(v_{n}t_{n}v_{n}u_{n}),\ \ x=p(w_{k,n}r_{k,n}w_{k,n}u_{ k,n}).$$

 There is a composition of disc exchanges contributing to both $\psi (\gamma _{1})$ and $\psi (\zeta (y))$. We consider the disc exchanges affecting $\gamma _{1}$, which are similar to those for \ref{2.7}. However, there is no exchange for the word $a$ preceding $v_{n}u_{n}$ in $v_{n}t_{n}v_{n}u_{n}$.  As in the case $k=n=0$, there is  a difference from the proof of \ref{2.7}, because this time the disc $v_{k}D_{1}'$ contains the point $p(v_{n}u_{n})$. But this is the only point in $O(y)$ that $v_{k}D_{1}'$ contains, and it is not in $P(y)$. So we have a similar picture to Figure 6, and (\ref{2.8.1}) holds. 

For $\zeta (y)$, as before, the first disc exchange which intersects $\zeta (y)$ does not move the endpoint, but does create a hook round $R_{3}D'$ exactly as in the picture for the case $k=n=0$.  The next exchange --- the inverse image of the first one affecting $\gamma _{1}$ --- is along the lower unit circle and reduces the size of this hook, replacing it by one of the same shape. As in \ref{2.7},  and as for $\beta (y)$, the disc exchanges which affect $\psi _{\ell }(\zeta (y))$  then alternate between a group of exchanges on the upper and lower unit circle. Finally, for  $m=26.2^{n}-11$, we have $\psi _{m}(\zeta (y))=\zeta (x)$ up to homotopy preserving $P(x)$ and $\psi _{m}(\zeta (y))=\psi (\zeta (y))$ up to homotopy preserving $P(x)$. So (\ref{2.8.2}) holds as before.

\Box

\section{Remarks}\label{3}
   
The pair $(a,b)=(L_{3}(L_{2}R_{3})^{2},L_{3}^{5})$ can be varied, to 
$$(L_{3}(L_{2}R_{3}) ^{2j},L_{3}(L_{2}R_{3})^{2j-2}L_{3}^{4})$$
 for any $j\geq 1$ or 
 $$(L_{3}(L_{2}R_{3})^{2j-1}L_{3}^{2}L_{2}R_{3},L_{3}(L_{2}R_{3})^{2j}L_{3}^{2})$$
 for any $j\geq 1$.
 
 Let us call such an ordered  pair $(a,b)$ an {\em{exchangeable pair}}. Note that each $a$ and each $b$ has an odd number of letters $L_{3}$ and $L_{2}$. Then if we examine the proof of \ref{2.7}, which works more generally than \ref{2.1}, we see that  if  $y=p(eau)$  then $\sigma _{\beta (y)}\circ s$
 and $\sigma _{\beta (x)}\circ s$ are Thurston equivalent  for some $x\in D(eb)$, whenever the following conditions hold.
 \begin{itemize}
\item[1.] $D(e)\subset D(L_{3}^{2\ell +1}L_{2})$ for some $\ell \geq 0$.
 \item[2.] $e$ consists only of the letters $L_{3}$, $L_{2}$ and $R_{3}$, with an even number of letters $L_{3}$ and $L_{2}$.
 \item[3.] $u$ consists only of the letters $L_{3}$, $L_{2}$ and $R_{3}$, apart from the last letter $C$.
  \item[4.] No suffix $u'$ of $eau$ satisfies $D(u')\subset D(eb)$.
 \item[5.] $eau<u<eb$, but no suffix $u'$ of $eau$ with $\vert u'\vert >\vert u\vert $ satisfies $eau<u'<eb$ .
 \item[6.] Any suffix $u'$ of $u$ with $\vert u'\vert <\vert u\vert $, which satisfies $eaebu<u'<ebu$, satisfies the stronger condition
 $eau<u'<ebu$.
 \item[7.] For any suffix $u'$ of $u$ which satisfies $eau<u'<eb$, the prefix of $u'$ in $u$ is $a_{1}$, for some exchangeable pair $(a_{1},b_{1})$.
  \item[8.] If $e'a'$ is a proper suffix of $e$ for some exchangeable pair $(a',b')$, such that $D(eau)$ is between $D(e'a'eau)$ and $D(e'b'eau)$, then $e'$ has an even number of letters $L_{3}$ and $L_{2}$. Note that this then implies that $e'a'eau<eau<e'b'eau$.

 \end{itemize}
 
 Once these conditions are satisfied, we can define $u_{1}$ to be obtained from $u$ by replacing $a_{1}$ by the second coordinate $b_{1}$ of the exchangeable pair, for every suffix $u'$ of $u$ with $eau<u<eb$, where $a_{1}$ is the first coordinate of the exchangeable pair which precedes $u'$. Then $\sigma _{\beta (y)}\circ s$
 and $\sigma _{\beta (x)}\circ s$ are Thurston equivalent, where $x=p(ebu_{1})$. The idea in \ref{2.7} is to have a word $v$ which can be written in the form $eau$ in many different ways.
 
This clearly gives a lot of choice, but it is unsatisfactory that the conditions imposed on $u$ are only satisfied by a zero density set of words, as the length tends to infinity. One would like to be able to say that, for a fixed choice of $e$, and for $u$ in a set of postive density, there exists $x\in D(eb)$ such that $\sigma _{\beta (y)}\circ s$
 and $\sigma _{\beta (x)}\circ s$ are Thurston equivalent. Such questions will be addressed in \cite{R6}. 
 
 I am not sure whether or not the growth rate of $n$ equivalent matings or captures with critical point of period or preperiod of $O(2^{n})$ is the best possible. A related question is on the growth of the length of words $v$ which can be written as $eau$ for $n$ different choices of $(e,a,u)$ where $(e,a,u)$ satisfies the conditions listed above. This growth seems to be $O(2^{n})$. Once again, it is hoped to address the issue properly  in \cite{R6}.

\end{document}